\documentclass{article}

\usepackage[english]{babel}
\usepackage{amssymb}
\usepackage{amsmath}
\usepackage{amsthm}

\newtheorem{thm}{Theorem}[section]
\newtheorem{defi}[thm]{Definition}
\newtheorem{lem}[thm]{Lemma}
\newtheorem{prop}[thm]{Proposition}
\newtheorem{cor}[thm]{Corollary}
\newtheorem{rmk}[thm]{Remark}
\newtheorem{ex}[thm]{Example}

\newcommand{\hh}{{\mathbb{H}}}
\newcommand{\cc}{{\mathbb{C}}}
\newcommand{\rr}{{\mathbb{R}}}
\newcommand{\nn}{{\mathbb{N}}}
\newcommand{\zz}{{\mathbb{Z}}}
\newcommand{\s}{{\mathbb{S}}}
\newcommand{\sing}{{\mathcal{S}}}
\newcommand{\z}{{\mathcal{Z}}}
\newcommand{\dom}{{\mathcal{D}}}
\newcommand{\lef}{{\mathcal{L}}}
\newcommand{\rig}{{\mathcal{R}}}
\newcommand{\minuszero}{{\setminus \{0\}}}

\newcommand{\minusr}{{\setminus \{r\}}}

\title{\bf Poles of regular quaternionic functions}
\author{ Caterina Stoppato \\
\normalsize Dipartimento di Matematica ``Ulisse Dini'', Universit\`a di Firenze \\ 
\normalsize Viale Morgagni 67/A, 50134 Firenze, Italy\\
\normalsize stoppato@math.unifi.it\\}

\date{  }

\begin{document}
\maketitle


\begin{abstract}
This paper studies the singularities of Cullen-regular functions of one quaternionic variable, as defined in \cite{advances}. The quaternionic Laurent series prove to be Cullen-regular. The singularities of Cullen-regular functions are thus classified as removable, essential or poles. The quaternionic analogues of meromorphic complex functions, called semiregular functions, turn out to be quotients of Cullen-regular functions with respect to an appropriate division operation. This allows a detailed study of the poles and their distribution.
\end{abstract}


\section{Introduction}

Denote by $\hh$ the skew field of real quaternions. Recall that it is obtained by endowing $\rr^4$ with the multiplication operation defined on the standard basis $1,i,j,k$ by $$i^2 = j^2 = k^2 = -1,$$ $$ ij = -ji = k, jk = -kj = i, ki = -ik = j,$$ $$1^2 = 1, 1 i = i 1 = i, 1 j = j 1 = j, 1 k = k 1 = k$$ and extended by distributivity to all quaternions $q = x_0 + x_1 i + x_2 j + x_3 k$.  A new theory of quaternion-valued functions of one quaternionic variable has been proposed by G. Gentili and D. C. Struppa in \cite{cras, advances}. The theory is based on a definition of regularity for quaternionic functions inspired by C. G. Cullen \cite{cullen}. Several interesting results are proven in \cite{advances}, including the Cullen-regularity of quaternionic power series and some basic properties of their zeros. The study of the zero-sets has been further developed in \cite{zeros, multiplicity}.

Let us quickly review the definition of Cullen-regular function and the basic properties of such a function. Denote by $\s$ the two-dimensional sphere of quaternionic imaginary units: $\s = \{q \in \hh : q^2 =-1\}$. For all imaginary unit $I \in \s$, let $L_I = \rr + I \rr$ be the complex line through $0, 1$ and $I$.

\begin{defi}\label{definition}
Let $\Omega$ be a domain in $\hh$ and let $f : \Omega \to \hh$ be a real differentiable function. $f$ is said to be \textnormal{Cullen-regular} if, for all $I \in \s$, the function $\bar \partial_I f : \Omega \cap L_I \to \hh$ defined by
\begin{equation}
\bar \partial_I f (x+Iy) = \frac{1}{2} ( \frac{\partial}{\partial x}+I\frac{\partial}{\partial y} ) f_I (x+Iy)
\end{equation}
vanishes identically.
\end{defi}

With the notations $\Omega_I = \Omega \cap L_I$ and $f_I = f_{|_{\Omega_I}}$, we may refer to the vanishing of $\bar \partial_I f$ saying that the restriction $f_I$ is holomorphic on $\Omega_I$. The following result clarifies the meaning of such a condition.

\begin{lem} \label{splitting}
Let $f : \Omega \to \hh$ be a Cullen-regular function and choose imaginary units $I,J \in \s$ with $I \perp J$. There exist functions $F,G : \Omega_I \to L_I$ such that $f_I = F + G J$ on $\Omega_I$. With the natural identification between $L_I$ and the complex field $\cc$, the functions $F,G$ are holomorphic, complex-valued functions of one complex variable.
\end{lem}

From now on we will omit Cullen's name and refer to these functions just as regular functions. As observed in \cite{advances}, a quaternionic power series $\sum_{n \in \nn} q^n a_n$ with $a_n \in \hh$ defines a regular function in its domain of convergence, which proves to be a ball $B(0,R) = \{q \in \hh : |q| <R\}$. In the same paper, it is proven that

\begin{thm} 
If $f : B = B(0,R) \to \hh$ is regular then there exist quaternions $a_n \in \hh$ such that $f(q)=\sum_{n \in \nn} q^n a_n$ for all $q \in B$. In particular, $f \in C^{\infty}(B)$.
\end{thm}

We may thus identify the set of regular functions on $B(0,R)$ with the set $\dom_R$ of quaternionic power series converging in $B(0,R)$. In \cite{advances} many basic results in complex analysis are extended to regular functions of this type: the identity principle, the maximum modulus principle, the Cauchy representation formula, the Liouville theorem, the Morera theorem and the Schwarz lemma. A version of the open mapping theorem has been recently proven in \cite{openarxiv, open}. The peculiar properties of the zeros proven in \cite{advances, zeros}, which we summarize in section \ref{preliminarysection}, arouse a new question. Do these functions have point singularities resembling the poles of holomorphic complex functions? In section \ref{laurentsection}, we give a positive answer to this question: we prove that a quaternionic Laurent series $$\sum_{n \in \zz} q^n a_n$$ defines a regular function on its domain of convergence, which is a four-dimensional spherical shell $A(0, R_1, R_2) = \{q \in \hh : R_1 < |q| < R_2\}$. This allows us to construct functions which are regular on a punctured ball $B(0,R) \minuszero$ and have a singularity at $0$. Moreover, we prove that any function which is regular on a spherical shell $A(0, R_1, R_2)$ admits a Laurent series expansion centered at $0$. This result is extended in section \ref{typesection} to the following:

\begin{thm}\label{introexpansion}
Let $f$ be a regular function on a domain $\Omega$, let $p \in \hh$ and let $L_I$ be a complex line through $p$. If $\Omega$ contains an annulus $A_I = A(p,R_1,R_2) \cap L_I$ then there exist $\{a_n\}_{n \in \zz} \subseteq  \hh$ such that $f_I(z) = \sum_{n \in \zz} (z-p)^n a_n$ for all $z \in A_I$. If, moreover, $p \in \rr$ then $f(q) = \sum_{n \in \zz} (q-p)^n a_n$ for all $q \in A(p,R_1,R_2) \cap \Omega$.
\end{thm}

This allows us to define quaternionic analogues of the concepts of \textit{pole} and  \textit{essential singularity} of holomorphic functions of one complex variable. We define the (classical) order $ord_f(p)$ of a pole and call a function $f$  \textit{semiregular} if it does not have essential singularities or, equivalently, if the restriction $f_I$ is meromorphic for all $I \in \s$. Note that, by the final statement of theorem \ref{introexpansion}, real singularities are completely analogous to singularities of holomorphic functions of one complex variable. There is no resemblance to the case of several complex variables: no such result as Hartog's lemma can hold. As for non-real singularities, we remark that theorem \ref{introexpansion} only provides information on the complex line $L_I$ through the point $p$; we apparently cannot predict the behavior of the function in a (four-dimensional) neighborhood of $p$. In order to overcome this difficulty, in section \ref{divisionsection} we introduce some new algebraic manipulation. We associate to any couple of regular functions $g,h: B(0,R) \to \hh$ a function $h^{-*}*g$, called the  \textit{left regular quotient} of $h$ and $g$. We study the basic properties of such a function and prove that it is semiregular on $B(0,R)$. We are thus able to conclude that if $f_I(z) = (z-p)^{-1} g_I(z)$ for $g_I$ holomorphic on a disk $B_I(0,R) = B(0,R) \cap L_I$ containing $p$ then $$f(q) = (q-p)^{-*}*g(q)$$ for some function $g$ which is regular on $B(0,R)$. As a consequence, in section \ref{semiregularsection} we prove the following result.

\begin{thm}
Let $f$ be a semiregular function on $B(0,R_0)$. For all $R<R_0$, $I \in \s$ there exist a polynomial $P(q)$ having coefficients in $L_I$ and a regular function $g : B(0,R) \to \hh$ such that $f = P^{-*}*g$ on $B(0,R)$.
\end{thm}

In particular $f$ is semiregular on $B(0,R_0)$ if and only if  $f_{|_{B(0,R)}}$ is a left regular quotient for all $R<R_0$. This allows the definition of a multiplication operation $*$ on the set of semiregular functions on a ball and the proof of the following result (where we denote $h^{*n} = h*...*h = *_{j=1}^n h$ the nth power of a regular function $h$ with respect to $*$-multiplication).

\begin{thm}
Let $f$ be a semiregular function on $B = B(0,R)$, choose $p \in B$ and let $m = ord_f(p), n = ord_f(\bar p)$. There exists a unique semiregular function $g$ on $B$ such that
\begin{equation}
f(q) = [(q-p)^{*m}*(q-\bar p)^{*n}]^{-*}* g(q)
\end{equation}
The function $g$ is regular near $p$ and $\bar p$ and $g(p) \neq 0, g(\bar p) \neq 0$, provided $m>0$ or $n>0$.
\end{thm}

The previous result allows the study of the structure of the poles:

\begin{thm}[Structure of the poles]
If $f$ is a semiregular function on $B = B(0,R)$ then $f$ extends to a regular function on $B$ minus a union of isolated real points or isolated 2-spheres of the type $x+y\s = \{x+yI : I \in \s\}$ with $x,y \in \rr, y \neq 0$. All the poles on each 2-sphere $x+y\s$ have the same order with the possible exception of one, which must have lesser order.
\end{thm}

Finally, in section \ref{multiplicitysection} we prove the following.

\begin{thm}
Let $f$ be semiregular on $B = B(0,R)$ and suppose $f \not \equiv 0$. For all $x+y\s \subseteq B$, there exist $m \in \zz, n \in \nn$, $p_1,...,p_n \in x+y\s$ with $p_i \neq \bar p_{i+1}$ for all $i,j$ such that
\begin{equation}
f(q) = [(q-x)^2+y^2]^m (q-p_1)*(q-p_2)*...*(q-p_n)*g(q)
\end{equation}
for some semiregular function $g$ on $B$ which does not have poles nor zeros in $x+y\s$.
\end{thm}

This theorem allows to extend to transcendental functions the concepts of \textit{spherical multiplicity} and \textit{isolated multiplicity} of the zeros defined in \cite{multiplicity} for polynomials. It also leads to analogous definitions for the poles of a semiregular function.


\section{Preliminary results}\label{preliminarysection}

We now run through the basic properties of the zero-sets of regular functions. In \cite{advances} it is proven that

\begin{thm}\label{symmetry}
Let $f : B(0,R) \to \hh$ be a regular function and let $x, y \in \rr$ be such that $x^2 +y^2 <R^2$. If there exist distinct imaginary units $I,J \in \s$ such that $f(x+yI) = f(x+yJ) = 0$, then $f(x+yK) = 0$ for all $K \in \s$.
\end{thm}

In other words, if $f$ has more than one zero on $x+y\s = \{x+yI : I \in \s\}$ then it vanishes identically on $x+y\s$. Note that $x+y\s$ is a 2-sphere if $y \neq 0$, a real singleton $\{x\}$ if $y= 0$.

\begin{ex}
The polynomial function $f(q) = q^2+1$ vanishes on $\s$. For all $x,y \in \rr$, the function $g(q) = (q-x)^2 +y^2 = q^2 - q 2x + x^2+y^2$ vanishes on $x+y\s$.
\end{ex}

In \cite{zeros} the zero-set is further characterized as follows.

\begin{thm}\label{structure}
Let $f$ be a regular function on an open ball $B(0,R)$. If $f$ is not identically zero then its zero-set consists of isolated points or isolated 2-spheres of the form $x + y \s$, for $x,y \in \rr, y\neq 0$.
\end{thm}

This result and the previous are proven for polynomials in \cite{shapiro} using quite simple tools. On the contrary, the study of the zero-set conducted in \cite{zeros} requires the introduction of the following operations on regular functions $f : B(0,R) \to \hh$.

\begin{defi}\label{multiplication}
Let $f, g$ be regular functions on an open ball $B = B(0,R)$ and let $f(q) = \sum_{n \in \nn} q^n a_n, g(q) = \sum_{n \in \nn} q^n b_n$ be their power series expansions. We define the \textnormal{regular product} of $f$ and $g$ as the regular function $f*g : B \to \hh$ defined by
\begin{equation}
f*g(q) = \sum_{n \in \nn} q^n c_n, \ c_n = \sum_{k=0}^n a_k b_{n-k}.
\end{equation}
Moreover, we define the \textnormal{regular conjugate} of $f$, $f^c : B \to \hh$, by $f^c(q) = \sum_{n \in \nn} q^n \bar a_n$ and the \textnormal{symmetrization} of $f$, as $f^s = f * f^c = f^c*f$. 
\end{defi}

The series $f*g$ and $f^c$ clearly converge in $B$. Also note that $f^s(q) = \sum_{n \in \nn} q^n r_n$ with $r_n = \sum_{k = 0}^n a_k \bar a_{n-k} \in \rr$. Since no confusion can arise, we may also write $f(q)*g(q)$ for $f*g(q)$.

\begin{rmk}
Fix $R$ with $0<R\leq \infty$ and let $\dom_R$ be the set of regular functions $f:B(0,R)\to \hh$. Then $(\dom_R,+,*)$ is an associative algebra over $\rr$. 
\end{rmk}

As observed in \cite{zeros}, the zeros of regular functions cannot be factored with respect to the standard multiplication of $\hh$. However, a factorization property is proven in \cite{zeros} in terms of $*$-multiplication, extending the results proven for polynomials in \cite{lam, serodio}. 

\begin{thm}\label{factorization}
Let $f : B = B(0,R) \to \hh$ be a regular function and let $p \in B$. Then $f(p) = 0$ if and only if there exists another regular function $g : B \to \hh$ such that $f(q) = (q-p) * g(q)$.
\end{thm}

\begin{ex}
Fix $x,y \in \rr$ and consider $f(q) = (q-x)^2 +y^2$. For all $I \in \s$, $f$ vanishes at $p = x+yI$ and it can be factored as $f(q) = q^2 - q 2x + x^2+y^2 = q^2 - q (p + \bar p) + p \bar p = (q-p)*(q-\bar p)$.
\end{ex}

If we define the nth regular power of $f$ as $f^{*n} = f*...*f = *_{i=1}^n f$, we can give the following definition.

\begin{defi}\label{multiplicity}
Let $f : B = B(0,R) \to \hh$ be a regular function, suppose $f \not \equiv 0$ and let $p \in B$. We define the \textnormal{(classical) multiplicity} of $p$ as a zero of $f$ and denote by $m_f(p)$ the largest $n \in \nn$ such that there exists a regular function $g : B \to \hh$ with $f(q)=(q-p)^{*n}*g(q)$.
\end{defi}

The classical multiplicity is a consistent generalization of the complex multiplicity; in other words, if $p \in L_I$ then $m_f(p)$ is the largest $n \in \nn$ such that there exists a holomorphic function $g_I$ with $f_I(z) = (z-p)^{n}g_I(z)$. This proves, a posteriori, that definition \ref{multiplicity} is well posed: indeed, by the identity principle proven in \cite{advances}, $f_I \not \equiv 0$; thus we cannot factor $z-p$ out of $f_I(z)$ ``infinitely many times''. Finally, the zero set of a regular product is completely characterized in terms of the zeros of the two factors by the following result, which is proven in \cite{zeros} and extends \cite{lam, serodio}.

\begin{thm}\label{zerosmultiplication}
Let $f,g$ be regular functions on an open ball $ B = B(0,R)$ and let $p \in B$. If $f(p) = 0$ then $f*g(p) = 0$, otherwise $f*g(p) = f(p)(f(p)^{-1} p f(p))$. In particular $p$ is a zero of $f*g$ if and only if $f(p) = 0$ or $g(f(p)^{-1} p f(p))=0$.
\end{thm}

We conclude this section recalling that the zero-sets of $f^c$ and $f^s$ are characterized in \cite{zeros} as follows.

\begin{thm}\label{conjugatezeros}
Let $f$ be a regular function on $B = B(0,R)$. For all $x,y \in \rr$ with $x+y\s \subseteq B$, the zeros of the regular conjugate $f^c$ on $x+y\s$ are in one-to-one correspondence with those of $f$. Moreover, the symmetrization $f^s$ vanishes exactly on the sets $x+y\s$ on which $f$ has a zero.
\end{thm}


\section{Laurent series and expansion}\label{laurentsection}

The first step in the study of point singularities is generalizing the theory of Laurent series to the quaternionic case. The domain of convergence of a quaternionic Laurent series $\sum_{n \in \zz} q^n a_n$ is a four-dimensional spherical shell $A(0, R_1, R_2) = \{q \in \hh : R_1 < |q| < R_2\}$. More precisely one can prove, just as in the complex case, the following result.

\begin{lem}\label{convergence}
Let $\{a_n\}_{n \in \zz} \subseteq \hh$. There exist $R_1, R_2$ with $0 \leq R_i \leq \infty$ such that
 \begin{enumerate}
\item the series  $\sum_{n \in \nn} q^n a_n$ and $ \sum_{n \in \nn} q^{-n} a_{-n}$ both converge absolutely and uniformly on compact subsets of $A = A(0,R_1,R_2)$;
\item for all $q \in \hh \setminus \bar A$ (except possibly $0$ if $A = \emptyset$), either $\sum_{n \in \nn} q^n a_n$ or $ \sum_{n \in \nn} q^{-n} a_{-n}$ diverge.
\end{enumerate}
\end{lem}

Note that $A=A(0,R_1,R_2)$ is empty if and only if $R_1\geq R_2$. If this is not the case then we define the sum of the series $\sum_{n \in \zz} q^n a_n$ as $\sum_{n \in \nn} q^n a_n +  \sum_{n>0} q^{-n} a_{-n}$ for all $q$ in $A$, which we may call the domain of convergence of the series.

\begin{thm}\label{laurentregularity}
Let $\sum_{n \in \zz} q^n a_n$ have domain of convergence $A = A(0,R_1,R_2)$ with $R_1 < R_2$. Then $f : A \to \hh \ \ q \mapsto \sum_{n \in \zz} q^n a_n$ is a regular function.
\end{thm}

The proof follows by computation from definition \ref{definition}. We will now prove that all regular functions $f : A(0,R_1,R_2) \to \hh$ admit Laurent series expansions. In order to prove it, we make use of the identity principle. This result is proven in \cite{advances} for functions which are regular on a ball $B(0,R)$, but it easily extends to a larger class of domains.

\begin{prop}[Identity Principle]\label{identity}
Let $\Omega$ be a domain in $\hh$ intersecting the real axis and having connected intersection $\Omega_I = \Omega \cap L_I$ with any complex line $L_I$. If $f,g : \Omega \to \hh$ are regular functions which coincide on $\Omega \cap \rr$ then they coincide on the whole domain $\Omega$.
\end{prop}

\begin{proof} Let $h = f-g$ and let us prove $h \equiv 0$ on $\Omega$. Choose any imaginary unit $I \in \s$ and consider the restriction $h_I = h_{|_{\Omega_I}}$. Since $h_I : \Omega_I \to \hh$ is holomorphic and it vanishes on the set $\Omega \cap \rr$, which is not discrete, $h_I$ must vanish identically on $\Omega_I$.
\end{proof}

We are now ready to prove the following.

\begin{thm}[Laurent Series Expansion]\label{expansion}
Let  $A = A(0,R_1,R_2)$ with $0 \leq R_1< R_2$ and let $f : A \to \hh$ be a regular function. There exist $\{a_n\}_{n \in \zz} \subseteq  \hh$ such that
\begin{equation}
f(q) = \sum_{n \in \zz} q^n a_n
\end{equation}
for all $q \in A$.
\end{thm}

\begin{proof}
Choose a complex line $L_I$ and consider the annulus we get by intersecting $L_I$ with the shell $A$: $$A_I = A_I(0,R_1,R_2) = \{z \in L_I : R_1<|z| < R_2\}.$$ Consider the restriction $f_I = f_{|_{A_I}}$ and choose $J \in \s, J \perp I$. As we saw in lemma \ref{splitting}, we can find two functions $F, G : A_I \to L_I$ which are holomorphic (with the natural identification between $L_I$ and $\cc$) and such that $f_I = F + G J$. Let $F(z) = \sum_{n \in \zz} z^n \alpha_n$ and $G(z) =  \sum_{n \in \zz} z^n \beta_n$ be the Laurent series expansions of the functions $F$ and $G$ (which have coefficients $\alpha_n, \beta_n \in L_I$). If we let $a_n = \alpha_n + \beta_n J$ for all $n \in \zz$, then $$f_I(z) =  \sum_{n \in \zz} z^n a_n$$ for all $z \in A_I$. Now consider the quaternionic Laurent series $\sum_{n \in \zz} q^n a_n$.  By lemma \ref{convergence} it converges in $A$. Hence, by theorem \ref{laurentregularity}, its sum defines a regular function on $A$. This function coincides with $f$ on $A_I$ by construction. We can conclude, using the identity principle \ref{identity}, that it coincide with $f$ on the whole domain $A$.
\end{proof}

The above argument also proves that

\begin{thm}\label{existence}
Let  $A_I = A_I(0,R_1,R_2)$ with $I \in \s, 0 \leq R_1< R_2$ and let $f_I : A_I \to \hh$ be holomorphic. There exists exactly one regular function $g : A(0,R_1,R_2) \to \hh$ such that $f_I(z) = g_I(z)$ for all $z \in A_I$.
\end{thm}


\section{Types of singularities}\label{typesection}

We now remark that the results proven in the previous section not only work for point $0$. It is easy to prove that for all $r \in \rr$ the series $\sum_{n \in \zz} (q-r)^n a_n$ converges and defines a regular function on some shell $A(r,R_1,R_2) = \{q \in \hh : R_1<|q-r|<R_2 \}$. Thus we can generalize all the above formulae just by substituting $r$ to $0$, $(q-r)^n$ to $q^n$ et cetera. On the contrary, we cannot apply the same procedure to a non-real quaternion $p \in \hh \setminus \rr$. Indeed, a quaternionic Laurent series centered at such a point $p$ has a nice convergence domain, $A(p,R_1,R_2) = \{q \in \hh : R_1<|q-p|<R_2 \}$, but its sum does not define, in general,  a regular function. Indeed it is easy to check that (due to the non-commutativity of $\hh$) the function $P(q) = (q-p)^n$ is not regular for $n \in \zz \setminus \{0,1\}$ and $p \in \hh \setminus \rr$. Nevertheless, the first part of the proof of theorem \ref{expansion} still works at a non-real point. This leads to the following result:

\begin{thm}\label{weakexpansion}
Let $f$ be a regular function on a domain $\Omega$, let $p \in \hh$ and let $L_I$ be a complex line through $p$. If $\Omega$ contains an annulus $A_I = A(p,R_1,R_2) \cap L_I$ with $0 \leq R_1< R_2$ then there exist $\{a_n\}_{n \in \zz} \subseteq  \hh$ such that
\begin{equation}\label{weakexpansioneq}
f_I(z) = \sum_{n \in \zz} (z-p)^n a_n
\end{equation}
for all $z \in A_I$. If, moreover, $p =r \in \rr$ then $f(q) = \sum_{n \in \zz} (q-r)^n a_n$ for all $q \in A(r,R_1,R_2) \cap \Omega$.
\end{thm}

For all $p$ and all complex lines $L_I$ through $p$, this result allows us to classify the behavior of $f(q)$ when $q$ approaches $p$ along $L_I$. Note that if $p$ does not lie in $\rr$ then there exists exactly one complex line $L_I$ through $p$. If, on the contrary, $p$ is real then it belongs to all complex lines $L_I$; however, in this case the coefficients $a_n$ which appear in equation \ref{weakexpansioneq} are the same for all $I \in \s$, thanks to the final statement of theorem \ref{weakexpansion}. We can thus give the following definition.

\begin{defi}
Let $f$, $p$ and $a_n$ be as in theorem \ref{weakexpansion}. The point $p$ is called a \textnormal{pole} if there exists an $n \in \nn$ such that $a_{-m} = 0$ for all $m>n$; the minimum of such $n \in \nn$ is called the  \textnormal{(classical) order} of the pole and denoted as $ord_f(p)$. If $p$ is not a pole for $f$ then we call it an \textnormal{essential singularity} for $f$.
\end{defi}

When no confusion can arise, we omit the adjective ``classical'' for the sake of simplicity. For a real point we derive from theorem \ref{weakexpansion} the following classification.

\begin{thm}\label{realclassification}
Let  $A = A(r,R_1,R_2)$ with $r \in \rr, 0 \leq R_1< R_2$ and let $f : A \to \hh$ be a regular function. 
\begin{enumerate}
\item If $r$ is a pole of order $0$ then $f$ extends to $B(r,R_2) = \{q \in \hh : |q-r|<R_2\}$ as a regular function.
\item If $r$ is a pole of order $n>0$ then there exists a regular function $g : B = B(r,R_2) \to \hh$ such that $f(q) = (q-r)^{-n}g(q)$ for all $q \in A$. In particular, $f$ extends to $B \minusr$ as a regular function and $\lim_{q \to r} |f(q)| = + \infty$.
\item Suppose $r$ to be an essential singularity. If $f$ extends to $B(r,R_2) \minusr$, then the modulus $|f|$ is unbounded on $U \minusr$ for all neighborhood $U$ of $r$; moreover, the limit $\lim_{q \to r} f(q)$ is not defined.
\end{enumerate}
\end{thm}

Note that when $p$ is not real the classification of $p$ as a pole or an essential singularity only depends on the restriction $f_I$ to the complex line through $p$. A priori, it does not predict the behavior of $f$ in a four-dimensional neighborhood of $p$. For instance, if $p$ is a pole of order $m$ then we derive from theorem \ref{weakexpansion} the existence of a holomorphic function $g_I : A_I \to \hh$ such that $$f_I(z) = \frac {1}{(z-p)^{m}} g_I(z).$$ However, for $p \in \hh \setminus \rr$ we cannot conclude that $f(q)$ equal $(q-p)^{-m} g(q)$ for some regular $g$ (note that the second expression does not generally define a regular function). In order to prove a result of this type we need an adequate division operation on regular functions, which we will define in section \ref{divisionsection}. Remark that even when $m=0$, i.e. $f_I$ extends to a holomorphic function on a neighborhood of $p$ in $L_I$, we may not conclude that $f$ extends to a regular function on a neighborhood of $p$ in $\hh$: theorem \ref{existence} was only proven for spherical shells $A(0,R_1,R_2)$ centered at $0$ and it does not immediately generalize to $A(p,R_1,R_2)$ with $p \in \hh \setminus \rr$. This is why we did not call a pole of order $0$ a removable singularity. Our caution will prove correct in section \ref{semiregularsection}. We will instead give the following, natural definition.

\begin{defi}
The point $p$ is called a \textnormal{removable singularity} if $f$ extends to a neighborhood of $p$ in $\hh$ as a regular function.
\end{defi}

We conclude this section defining an analogue to the concept of meromorphic function. We will call semiregular a function which does not have essential singularities. More precisely:
\begin{defi}
Let $\Omega$ be a domain in $\hh$, let $\sing \subseteq \Omega$ and suppose the intersection $\sing_I = \sing \cap L_I$ to be a discrete subset of $\Omega_I = \Omega \cap L_I$ for all $I \in \s$. A regular function $f : \Omega \setminus \sing \to \hh$ is said to be \textnormal{semiregular} on $\Omega$ if $\sing$ does not contain essential singularities for $f$.
\end{defi}

In other words, $f$ is semiregular on $\Omega$ if and only if, for all $I \in \s$, the restriction $f_I$ is a meromorphic function on $\Omega_I$. Note that we did not ask for the set of singularities $\sing$ to be discrete: in order to classify a point $p \in \sing$ as a pole or an essential singularity for $f$ it is enough for $p$ to be isolated in $\sing_I = \sing \cap L_I$ (see the hypotheses of theorem \ref{weakexpansion}). In section \ref{semiregularsection} we will present a detailed study of the functions which are semiregular on a ball $B(0,R)$. This study requires the introduction of the above mentioned division operation, which we undertake in the next section.


\section{Regular quotients}\label{divisionsection}

As we saw in section \ref{preliminarysection}, the zeros of a regular function cannot be factored with respect to the multiplication of $\hh$, but they have nice multiplicative properties in terms of the non-standard multiplication $*$. Similarly, the apparent difficulties we found in dealing with non-real poles can be solved in terms of a non-standard division operation. Denote by  $\z_f = \{q \in B(0,R) : f(q) = 0\}$ the zero-set of a function $f : B(0,R) \to \hh$. 

\begin{defi}
Let $f,g : B = B(0,R) \to \hh$ be regular functions and let $f^c,f^s$ be the regular conjugate and the symmetrization of $f$. We define the \textnormal{left regular quotient} of $f$ and $g$ as the function $f^{-*} * g : B \setminus \z_{f^s} \to \hh$ such that
\begin{equation}
f^{-*} * g (q) = \frac{1}{f^s(q)} f^c * g(q).
\end{equation}
The  \textnormal{right regular quotient} of $g$ and $f$ is the function $g*f^{-*} : B \setminus \z_{f^s} \to \hh$ defined by $g*f^{-*} =f^{-s}(g *f^c)$. Finally, we define the \textnormal{regular reciprocal} of $f$ as the function $f^{-*} = f^{-*} * 1 = 1*f^{-*}$.
\end{defi}

Since no confusion can arise, we will often write $(f(q))^{-*}$ for $f^{-*}(q)$. We will also use the shorthand notation $f^{-s}(q)$ for $\frac{1}{f^s(q)}$. Regular quotients are regular on their domains of definition by the following lemma, which can be proven by direct computation.

\begin{lem}
Let $f,g: B = B(0,R) \to \hh$ be regular functions and suppose the power series expansions of $g$ at $0$, $g(q) = \sum_{n \in \nn} q^n r_n$, has real coefficients $r_n \in \rr$. Then the function $h: B \setminus \z_g \to \hh$ defined by $h(q) = \frac{1}{g(q)} f(q)$ is regular.
\end{lem}

Moreover, left and right regular quotients of regular functions on a ball $B(0,R)$ are semiregular on $B(0,R)$:

\begin{thm}
Let $f,g : B = B(0,R) \to \hh$ be regular functions and consider the left quotient $f^{-*}*g : B \setminus \z_{f^s} \to \hh$. Each $p \in \z_{f^s}$ is a pole of order 
\begin{equation}
ord_{f^{-*}*g}(p) \leq m_{f^s}(p)
\end{equation}
for $f^{-*}*g$, where $m_{f^s}(p)$ denotes the classical multiplicity of $f^s$ at $p$. The same holds for the right quotient $g * f^{-*}$. In particular $f^{-*}*g$ and $g * f^{-*}$ are semiregular on $B$.
\end{thm}

\begin{proof}
Let $p = x+yI \in L_I$. If $m_{f^s}(p) = n$ then, as a consequence of theorem \ref{conjugatezeros}, $m_{f^s}(\bar p) = n$; there exists a holomorphic function $h_I$ with $h_I(p) \neq 0$ such that $f_I^s(z) = (z-p)^n(z-\bar p)^n h_I(z)=  \left[(z-x)^2+y^2\right]^n h_I(z)$. We observe that, since $f^s$ is a series with real coefficients, $h_I(z)$ must be have real coefficients, too. As a consequence, $f_I^s(z) = h_I(z) \left[(z-x)^2+y^2\right]^n = h_I(z) (z-\bar p)^n (z-p)^n$ and $$(f^{-*}*g)_I(z) = f_I^{-s}(z) (f^c *g)_I(z) = (z-p)^{-n}(z-\bar p)^{-n} h_I(z)^{-1}(f^c *g)_I(z)$$ where $(z-\bar p)^{-n} h_I(z)^{-1}(f^c *g)_I(z)$ is holomorphic in a neighborhood of $p$ in $L_I$. Moreover $(g*f^{-*})_I(z) = (z-p)^{-n}(z-\bar p)^{-n} h_I(z)^{-1}(g*f^c)_I(z)$, with $(z-\bar p)^{-n} h_I(z)^{-1}(g*f^c)_I(z)$ holomorphic in a neighborhood of $p$ in $L_I$.
\end{proof}

The regular quotient $f^{-*} * g(q)$ is related to the quotient $f(q)^{-1} g (q)= \frac{1}{f(q)} g(g)$ by the following result. First remark that, as a consequence of theorem \ref{conjugatezeros}, $\z_{f^c} \subseteq \z_{f^s}$, so that $f^c(q) \neq 0$ for all $q \in B \setminus \z_{f^s}$.

\begin{thm}\label{quotients}
Let $f,g$ be regular functions on $B=B(0,R)$. If we define $T_f : B \setminus \z_{f^s} \to \hh$ as $T_f(q) = f^c(q)^{-1} q f^c(q)$, then
\begin{equation}
f^{-*}*g(q)= \frac{1}{f( T_f(q))} g (T_f(q))
\end{equation}
for all $q \in B \setminus \z_{f^s}$.
\end{thm}

\begin{proof}
By theorem \ref{zerosmultiplication}, $f^c(q) * g(q) = f^c(q) g(T_f(q))$ for all $q \in B \setminus \z_{f^s}$. We conclude by computation: $$f^{-*}*g(q) = f^{-s}(q) f^c*g(q)= [f^c*f(q)]^{-1} f^c*g(q)=$$ $$= [f^c(q) f(T_f(q))]^{-1} f^c(q)g(T_f(q)) = f(T_f(q))^{-1} f^c(q)^{-1} f^c(q) g(T_f(q)) = $$ $$= f(T_f(q))^{-1}g(T_f(q)).$$
\end{proof}

Since for all $x,y \in \rr, I \in \s$ and $p \in \hh \minuszero$ we have $p^{-1}(x+yI) p = p^{-1}x p + p^{-1} y I p = x+y J$ with $J = p^{-1} I p \in \s$, we remark that:

\begin{rmk}
For all $x,y \in \rr$ with $x^2+y^2<R^2$, the function $T_f$ maps the 2-sphere (or real singleton) $x+y\s = \{x+yI : I \in \s\}$ to itself.
\end{rmk}

In particular, since by theorem \ref{conjugatezeros} we have
$$\z_{f^s} = \bigcup_{x+yI \in \z_f} (x+y\s),$$
we conclude $T_f(B \setminus \z_{f^s}) \subseteq B \setminus \z_{f^s}$. Moreover, the following is proven in \cite{openarxiv, open}.

\begin{prop}
Let $f: B=B(0,R) \to \hh$ be a regular function. $T_f$ and $T_{f^c}$ are mutual inverses. In particular $T_f : B \setminus \z_{f^s} \to B \setminus \z_{f^s}$ is a diffeomorphism.
\end{prop}

Let us give an example.

\begin{ex}
For any fixed quaternion $p = x+yI \in \hh$, the regular reciprocal of the polynomial $f(q) = q-p$ is $$(q-p)^{-*} =  \frac{1}{(q-\bar p)*(q-p)} (q-\bar p)= \frac{1}{q^2 - q (p + \bar p) + \bar p p} (q-\bar p) = \frac{1}{(q- x)^2 + y^2} (q - \bar p),$$ where the polynomial $(q- x)^2 + y^2$ vanishes exactly on $x+y\s$. Moreover, by theorem \ref{quotients}, $$(q-p)^{-*} = \frac{1}{T_f(q) -p} = \frac{1}{(q-\bar p)^{-1} q (q-\bar p) -p} $$
In particular, since $T_f(q) = (q-\bar p)^{-1} q (q-\bar p) = (q-\bar p)^{-1} (q-\bar p) q = q$ for $q$ in the complex line $L_I$ through $p$, the function $(q-p)^{-*}$ coincides with $(q-p)^{-1} = \frac{1}{q-p}$ on $L_I \setminus \{p,\bar p\}$.
\end{ex}

We conclude this section explaining the algebraic meaning of regular quotients. In the complex case, the set of quotients $\frac{F}{G}$ of holomorphic functions $F,G$ on a disc $\Delta$ becomes a field when endowed with the usual operations of addition and multiplication. More precisely, it is the field of quotients of the integral domain (i.e. the commutative ring with no zero divisors) obtained by endowing the set of holomorphic functions $F$ on $\Delta$ with the natural addition and multiplication. As explained in \cite{rowen} (see also \cite{cohn, lam2}), the concept of field of quotients of an integral domain can be generalized to the non-commutative case as follows. 

\begin{thm}
We define a \textnormal{left Ore domain} as a domain (a ring with no zero divisors) $(D,+,\cdot)$ such that $Da \cap Db \neq \{0\}$ for all $a,b \in D \minuszero$. If this is the case, then the set of formal quotients $L = \{a^{-1}b : a,b \in D\}$ can be endowed with operations $+, \cdot$ such that: 
\begin{itemize}
\item[(i)] $D$ is isomorphic to a subring  of $L$ (namely $ \{1^{-1}a : a \in D\}$); 
\item[(ii)] $L$ is a skew field, i.e. a ring where every non-zero element has a multiplicative inverse (namely, $(a^{-1}b)^{-1} = b^{-1}a$).
\end{itemize}
The ring $L$ is called the \textnormal{classical left ring of quotients} of $D$ and, up to isomorphism, it is the only ring having the properties (i) and (ii). 
\end{thm}
On a \textit{right Ore domain} $D$, defined by $aD \cap bD \neq \{0\}$ for all $a,b \in D \minuszero$, we can similarly construct the \textit{classical right ring of quotients}. If $D$ is both a left and a right Ore domain, then (by the uniqueness property) the two rings of quotients are isomorphic and we may speak of the \textit{classical ring of quotients} of $D$.

\begin{prop}
Fix $R$ with $0<R\leq \infty$. The associative real algebra $(\dom_R,+,*)$ of regular functions on $B(0,R)$ is a left Ore domain and a right Ore domain. If we endow the set of left regular quotients $\lef_R = \{f^{-*}*g : f,g \in \dom_R, f \not \equiv 0\}$ with the multiplication $*$ defined by 
\begin{equation}
(f^{-*}*g)*(h^{-*}*k) = f^{-s}h^{-s} f^c*g*h^c*k
\end{equation}
then $(\lef_R,+,*)$ is a division algebra over $\rr$ and it is the classical ring of quotients of $(\dom_R,+,*)$. The same holds for $\rig_R = \{g*f^{-*} : f,g \in \dom_R, f \not \equiv 0\}$ with the multiplication defined by $(g*f^{-*})*(k*h^{-*}) = f^{-s}h^{-s} g*f^c*k*h^c$.
\end{prop}

\begin{proof} 
The multiplication $*$ is well defined on $\lef_R$: $f^{-*}*g = \tilde f^{-*}* \tilde g$ if and only if there exist $l, \tilde l$ such that $l * f = \tilde l * \tilde f$, $l*g=\tilde l*\tilde g$ and in this case we get by direct computation that $f^{-s}h^{-s} f^c*g*h^c*k = \tilde f^{-s}h^{-s} \tilde f^c*\tilde g*h^c*k$; the same can be done for the second factor $h^{-*}*k$. Clearly, $(\lef_R,+,*)$ is an associative algebra over $\rr$. We remark that $f^{-*}*g$ has inverse element $g^{-*}*f$ with respect to $*$: $$(f^{-*}*g) * (g^{-*}*f) = f^{-s}g^{-s} f^c*g*g^c*f = f^{-s}g^{-s} f^c*g^s*f =$$ $$= f^{-s}g^{-s} g^s f^c*f = f^{-s} f^c*f = f^{-s} f^s = 1$$ and, switching $f$ and $g$, $(g^{-*}*f) *(f^{-*}*g) = 1$. Thus $(\lef_R,+,*)$ is a division algebra. The same holds for $(\rig_R,+,*)$.

The ring $\dom_R$ is a domain, since $f*g \equiv 0$ iff $f \equiv 0$ or $g\equiv 0$. Moreover, $\dom_R$ is a left Ore domain: if $f,g \not \equiv 0$ then $(\dom_R*f) \cap (\dom_R*g)$ contains the non-zero element $f^s g^s = g^s f^s$, which can be obtained as $(g^s * f^c )*f$ or as $(f^s*g^c)*g$. Similarly,  $\dom_R$ is a right Ore domain. Thus the classical ring of quotients of $\dom_R$ is well defined. It must be isomorphic to both $\lef_R$ and $\rig_R$ by the uniqueness property: $\lef_R,\rig_R$ are skew fields which have $\dom_R$ as a subring and the inclusions $\dom_R \to \lef_R \ \ f \mapsto f = 1^{-*} *f$, $\dom_R \to \rig_R \ \ f \mapsto f = f*1^{-*}$ prove to be ring homomorphisms by direct computation.
\end{proof}


\section{Poles of semiregular functions}\label{semiregularsection}

We now prove that all functions which are semiregular in a neighborhood of a bounded ball $B$ can be expressed as left quotients on $B$.

\begin{thm}\label{polesfactorization}
Let $f$ be a semiregular function on $B(0,R_0)$ with $0<R_0\leq \infty$. Let $B = B(0,R)$ with $0<R<R_0$, choose $I \in \s$ and let $z_1,...,z_n$ be the poles of $f_I$ in $B_I = B \cap L_I$, listed according to their order $ord_f$. There exists a unique regular function $g : B \to \hh$ such that
\begin{equation}\label{polesfactorizationequation}
f(q) = \left[(q-z_1)*...*(q-z_n)\right]^{-*} * g(q)
\end{equation}
for all $q \in B$. Moreover, $g(z_j) \neq 0$ for all $j \in \{1,...,n\}$.
\end{thm}

\begin{proof}
As we observed in section \ref{typesection}, there exists a holomorphic $g_I : B_I \to \hh$, with $g_I(z_j) \neq 0$ for all $j$, such that $$f_I(z) = \frac{1}{(z-z_1)...(z-z_n)}g_I(z)$$ for all $z \in B_I \setminus\{z_1,...,z_n\}$. As a consequence of theorem \ref{existence}, $g_I$ extends to a regular $g : B \to \hh$. Consider the function $h(q) =  \left[(q-z_1)*...*(q-z_n)\right]^{-*} * g(q)$: it is regular on its domain of definition, which is $\Omega = B \setminus \bigcup_{j = 1}^n \left(x_j+y_j\s\right)$ if $z_j=x_j+y_jI_j$. We remark that, as a consequence of theorem \ref{quotients}, $h_I (z) =  [(z-z_1)...(z-z_n)]^{-1} g_I(z)=f_I(z)$ for all $z \in \Omega_I =  B_I \setminus\{z_1, \bar z_1,...,z_n, \bar z_n\}$. The identity principle \ref{identity} allows us to conclude that $f(q) = h(q)$ for all $q \in \Omega$.
\end{proof}

\begin{cor}
Let $0<R_0\leq \infty$. A function $f$ is semiregular on $B(0,R_0)$ if, and only if, $f_{|_{B(0,R)}} \in \lef_R$ for all $R<R_0$.
\end{cor}

We use the notation $f_{|_{B(0,R)}}$ for the sake of simplicity, instead of writing $f _{|_{\Omega \cap B(0,R)}}$ where $\Omega$ is the domain on which $f$ is regular.
Thanks to the previous corollary, we can define a multiplication operation $*$ on semiregular functions on $B(0,R_0)$ with $0<R_0\leq \infty$. Consider indeed two such functions $f,g$. On each ball $B(0,R)$ with $0<R<R_0$, the restrictions $f _{|_{B(0,R)}},g_{|_{B(0,R)}}$ can be represented as left regular quotients and we may consider their product $f _{|_{B(0,R)}}* g_{|_{B(0,R)}}$. Moreover, taking $R_2>R_1$ we get that $f _{|_{B(0,R_2)}}* g_{|_{B(0,R_2)}}$ equals $f_{|_{B(0,R_1)}} * g_{|_{B(0,R_1)}}$ on $B(0,R_1)$. We can thus define:

\begin{defi}
The \textnormal{(semi)regular product} of semiregular functions $f,g$ on $B(0,R_0)$ is the semiregular function $f*g$ on $B(0,R_0)$ such that $(f*g)_{|_{B(0,R)}} = f_{|_{B(0,R)}} * g_{|_{B(0,R)}}$ for all $R<R_0$. 
\end{defi}

We can now remark the following. Recall that we denote by $h^{*n} = h*...*h = *_{j=1}^n h$ the nth regular power of $h$. 

\begin{thm}\label{polefactorization}
Let $f$ be a semiregular function on $B(0,R_0)$ with $0<R_0\leq\infty$, choose $p = x+yI \in B(0,R_0)$ and let $m = ord_f(p), n = ord_f(\bar p)$. Without loss of generality, $m \leq n$. There exists a unique semiregular function $g$ on $B(0,R_0)$ such that
\begin{equation}
f(q) = [(q-p)^{*m}*(q-\bar p)^{*n}]^{-*}* g(q) = \left[(q-x)^2+y^2\right]^{-n} (q-p)^{*(n-m)}* g(q) 
\end{equation}
The function $g$ is regular in a neighborhood of $x+y\s$ and, if $n>0$, $g(p),g(\bar p) \neq 0$.
\end{thm}

\begin{proof}
For all $R<R_0$, the existence of a $g^{(R)} \in \lef_R$ such that $f(q) = [(q-p)^{*m}*(q-\bar p)^{*n}]^{-*}* g^{(R)} (q)$ on $B(0,R)$ is an immediate consequence of theorem \ref{polesfactorization}. Clearly, if $R_1<R_2$ then $g^{(R_1)}$ equals $g^{(R_2)}$ on $B(0,R_1)$. We can thus define a global $g$, semiregular on $B(0,R_0)$, such that $f(q) = [(q-p)^{*m}*(q-\bar p)^{*n}]^{-*}* g(q)$ on $B(0,R_0)$.
We conclude by observing that $$[(q-p)^{*m}*(q-\bar p)^{*n}]^{c} = (q-p)^{*n}*(q-\bar p)^{*m} =  [(q-x)^2+y^2]^{m}  (q-p)^{*(n-m)}$$ and $[(q-p)^{*m}*(q-\bar p)^{*n}]^{-s} = [(q-x)^2+y^2]^{-m-n}$, so that $$ [(q-p)^{*m}*(q-\bar p)^{*n}]^{-*}* g(q) = \left[(q-x)^2+y^2\right]^{-n} (q-p)^{*(n-m)}* g(q).$$
\end{proof}

We will soon use theorem \ref{polefactorization} to study the distribution of the poles. Let us first give two significant examples.

\begin{ex}
The regular function $f : \hh \setminus \s \to \hh$ defined by $$f(q) = (q^2+1)^{-*} = \frac{1}{q^2+1}$$ has a pole of order $1$ at any point $I \in \s$: we indeed have $f_I(z) =  \frac{1}{z-I}\frac{1}{z+I}$ for all $z \in L_I \setminus \{I,-I\}$.
\end{ex}

\begin{ex}
For any non-real quaternion $p = x+yI \in \hh \setminus \rr$, the function $f : \hh \setminus (x+y\s) \to \hh$ defined by $$f(q) = (q-p)^{-*}  = \left[(q-x)^2+y^2\right]^{-1} (q - \bar p)$$ has poles of order $1$ at all points of $x+y\s$ except $\bar p$, which has order $0$. Indeed, $f_I(z) = \frac{1}{(z- p)(z-\bar p)}(z-\bar p) =  \frac{1}{z -p}$ for all $z \in L_I  \setminus \{p,\bar p\}$, while for $p' = x+yJ$ with $J \in \s \setminus \{I, -I\}$ we have $f_J(z) = \frac{1}{(z- p')(z-\bar p')}(z-\bar p)$ for all $z \in L_J \setminus \{p',\bar p'\}$, where $p'-\bar p \neq 0$ and $\bar p' - \bar p \neq 0$.
\end{ex}

The previous example proves that 

\begin{rmk}
A pole of order $0$ is not always a removable singularity.
\end{rmk}

This is because $f(q) = (q-p)^{-*} $ does not extend as a regular function to a neighborhood of $\bar p$. Indeed any such neighborhood $U$ contains poles of order $1$ for $f$ and we conclude that $|f|$ is unbounded on $U$. We now study the distribution of the poles of a generic function which is semiregular on a ball.

\begin{thm}[Structure of the poles]\label{polestructure}
Let $f$ be a semiregular function on $B = B(0,R)$ with $0<R\leq\infty$. Then $f$ extends to a regular function on $B \setminus \sing$ with $\sing$ consisting of isolated 2-spheres (or real singletons) of the form $x+y\s$. All the poles on each sphere $x+y\s$ have the same order with the possible exception of one, which must have lesser order.
\end{thm}

\begin{proof}
Take $x+y\s \subseteq B$ and suppose that there exists $I \in \s$ such that $p = x+yI$ and $\bar p = x-yI$ have orders $m$ and $n$ with $m>0$ or $n>0$. By possibly substituting $-I$ to $I$, we may suppose $m \leq n$. By theorem \ref{polefactorization}, there exists a semiregular function on $B$ which is regular in a neighborhood $U$ of $x+y\s$ such that $f(q) = [(q-x)^2+y^2]^{-n} (q-p)^{*(n-m)}* g(q)$. Observing that the last expression is regular on $U \setminus (x+y\s)$ proves the first statement of the theorem.

We now prove the second statement. If we set $\tilde f(q) = (q-p)^{*(n-m)}* g(q)$ then $$f(q) = [(q-x)^2+y^2]^{-n} \tilde f(q),$$ $$f_J(z) = [z-(x+yJ)]^{-n} [z-(x-yJ)]^{-n} \tilde f_J(z)$$ for all $J \in \s$.
If $m<n$ then $\tilde f(x+yI) = 0$ and $\tilde f (x+yJ) \neq 0$ for all $J \in \s \setminus \{I\}$. The previous equation allows us to conclude $ord_f(x+yJ) = n$ for all $J \in \s \setminus \{I\}$. Since we know by hypothesis that $ord_f(x+yI)= ord_f(p) = m<n$, the thesis holds.
If $m = n$ then $\tilde f(x+yI) \neq 0$. If $\tilde f$ does not have zeros in $x +y \s$ then we conclude $ord_f(x+yJ) = n$ for all $J \in \s$. If, on the contrary, $\tilde f(x + yK)=0$ for some $K \in \s$ then we can factor $z-(x+yK)$ out of $\tilde f_K(z)$ and conclude that $ord_f(x + yK)<n$ while $ord_f(x + yJ)=n$ for all $J \in \s \setminus \{K\}$, as desired.
\end{proof}


\section{A different approach to multiplicity and order}\label{multiplicitysection}

The following peculiar property of the zeros of a quaternionic polynomial is shown in \cite{zeros}. Recall that $m_f(p)$ denotes the classical multiplicity of $p$ as a zero of a regular function $f$, defined in \ref{multiplicity} as the largest $n \in \nn$ such that $f(q) = (q-p)^{*n}*g(q)$ for some regular $g$.

\begin{prop}
Let $P(q)$ be a regular quaternionic polynomial of degree $d$ which does not have spherical zeros. Then $d \geq \sum_{q \in \z_P} m_P(q)$ and the inequality can be strict.
\end{prop}

If $P$ has a spherical zero $x+y \s$ then the situation is even more peculiar: clearly $m_P(q) > 0$ for all the (infinite) points $q \in x+y\s$. In order to overcome these apparent difficulties, an alternative approach has been recently introduced in \cite{multiplicity}. We may rephrase the definition given in \cite{multiplicity} as follows.

\begin{defi}\label{newmultiplicity}
Let $P(q)$ be a regular quaternionic polynomial and let $x,y \in \rr$. We say that $P$ has \textnormal{spherical multiplicity} $2m$ at $x+y\s$ if $m$ is the largest natural number such that $P(q) = [(q-x)^2+y^2]^m \tilde P(q)$ for some other polynomial $\tilde P$. If $\tilde P$ has a zero $p_1 \in x+y\s$ then we say that $P$ has \textnormal{isolated multiplicity} $n$ at $p_1$, where $n$ is the largest natural number such that there exist $p_2,...,p_n \in x+y\s$ and a polynomial $R(q)$ with $\tilde P(q) = (q-p_1)*(q-p_2)*...*(q-p_n)*R(q)$.
\end{defi}

As observed in \cite{multiplicity}, the previous definition yields:

\begin{prop}
If $P(q)$ is a regular quaternionic polynomial of degree $d$, then the sum of the spherical multiplicities and the isolated multiplicities of $P$ is $d$.
\end{prop}

The classical multiplicity $m_P$ of a polynomial $P$ is related to the spherical and isolated multiplicities of $P$ as follows:

\begin{rmk}\label{relationmult}
Let $P(q)$ be a regular quaternionic polynomial and let $p = x+yI \in \hh$. Then $P$ has spherical multiplicity $2 \min \{ m_P(p),m_P(\bar p)\}$ at $x+y\s$. Moreover, if $m_P(p)>m_P(\bar p)$, then $P$ has isolated multiplicity $n \geq m_P(p)-m_P(\bar p)$ at $p$.
\end{rmk}

Let us give some examples to clarify the previous remark.

\begin{ex}
The polynomial $P(q) = q^2+1$ vanishes on $\s$.  For all $I \in \s$, $P(q) = (q-I)*(q+I)$ has classical multiplicity $m_P(I) =1$ at $I$.  Moreover, $P$ has spherical multiplicity $2$ at $\s$.
\end{ex}

\begin{ex}
If $I \in \s$, then the polynomial $P(q) = (q-I) *(q-I) = (q-I)^{*2}$ only vanishes at $I$. $P$ has classical multiplicity $m_P(I) =2$ at $I$, it has spherical multiplicity $0$ at $\s$ and isolated multiplicity $2$ at $I$.
\end{ex}

\begin{ex}
If $I,J \in \s, I\neq J,I \neq -J$, then the polynomial $P(q) = (q-I) *(q-J) = q^2 - q (I+J) +IJ$ only vanishes at $I$, where it has classical multiplicity $m_P(I) = 1$  (see \cite{zeros} for details). Note that this is an example of polynomial having degree 2 greater that the sum $\sum_{q \in \z_P} m_P(q) = 1$ of the classical multiplicities of its zeros. According to definition \ref{newmultiplicity}, $P$ has spherical multiplicity $0$ at $\s$ and isolated multiplicity $2$ at $I$.
\end{ex}

We note that it is possible to combine the three cases presented above to build new examples. Fix $I \in \s$. We can easily construct a quaternionic polynomial $P$ of degree $d$ having classical multiplicities $m_P(I) = M, m_P(-I) = m$ at $I, -I$, spherical multiplicity $2m$ at $\s$ and isolated multiplicity $n$ at $I$ whenever $d,m,M,n \in \nn$ are such that $m \leq M$ and $M-m \leq n \leq d-m$. Thus remark \ref{relationmult} is sharp.

Definition \ref{newmultiplicity} does not immediately extend to transcendental functions: a priori, there may exist a regular function $f : B \to \hh$ and a 2-sphere $x+y\s \subseteq B$  such that we can factor out $q-p_j$ for ``infinitely many'' $p_j \in x+y\s$. We now prove that this is not the case. It is convenient to take care, at the same time, of the case of an $f$ having poles on $x+y\s$.

\begin{thm}\label{transfactorization}
Let $f$ be a semiregular function on $B = B(0,R)$, suppose $f \not \equiv 0$ and let $x+y\s \subseteq B$. There exist $m \in \zz, n \in \nn$, $p_1,...,p_n \in x+y\s$ with $p_i \neq \bar p_{i+1}$ for all $i$ such that
\begin{equation}\label{eqpoles}
f(q) = [(q-x)^2+y^2]^m (q-p_1)*(q-p_2)*...*(q-p_n)*g(q)
\end{equation}
for some semiregular function $g$ on $B$ which does not have poles nor zeros in $x+y\s$.
\end{thm}

\begin{proof} As we saw in the proof of theorem \ref{polefactorization}, if $x+y\s$ contains poles then there exists a $j> 0$ such that $f(q) = [(q-x)^2+y^2]^{-j} \tilde f(q)$ where $\tilde f$ does not have poles in $x+y\s$. Hence it is enough to prove the theorem for functions which are regular around $x+y\s$. This proof requires two steps:
\begin{enumerate}
\item[(i)] If $f$ is regular around $x+y\s$ and $f \not \equiv 0$ on $B$, then there exists an $m \in \nn$ such that $$f(q) = [(q-x)^2+y^2]^m h(q)$$ for some $h$ which does not vanish identically on $x+y\s$. Suppose indeed it were possible to find, for all $k \in \nn$, a function $h^{(k)}(q)$ such that $f(q) = [(q-x)^2+y^2]^k h^{(k)}(q)$. Then, choosing an $I \in \s$, the meromorphic function $f_I$ would have the factorization $f_I(z) = [(z-x)^2+y^2]^k h^{(k)}_I(z) = [z-(x+yI)]^k[z-(x-yI)]^k h^{(k)}_I(z)$ for all $k \in \nn$. This would imply $f_I \equiv 0$ and, by the identity principle \ref{identity}, $f \equiv 0$.
\item[(ii)] Let $h$ be a semiregular function on $B$ which does not have poles in $x+y\s$ nor vanishes identically on $x+y\s$. By theorem \ref{symmetry}, $g^{(0)}=h$ has at most one zero $p_1 \in x+y\s$. If this is the case then $h(q) = (q-p_1) * g^{(1)}(q)$  for some function $g^{(1)}$ which does not vanish identically on $x+y\s$. If for all $k \in \nn$ there existed a $p_{k+1} \in x+y\s$ and a $g^{(k+1)}$ such that $g^{(k)}(q) = (q-p_{k+1})*g^{(k+1)}$ then we would have $$h(q) = (q-p_1)*...*(q-p_k)*g^{(k)}(q)$$for all $k \in \nn$. This would imply, for the symmetrization $h^s$ of $h$, $$h^s(q) = [(q-x)^2+y^2]^k (g^{(k)})^s(q)$$ for all $k \in \nn$. By point (i), this would imply $h^s \equiv 0$. We could then conclude, applying theorems \ref{conjugatezeros} and \ref{structure}, that $h \equiv 0$, a contradiction. Thus there exists an $n \in \nn$ such that $g^{(n)}$ does not have zeroes in $x+y\s$ and, setting $g = g^{(n)}$, we conclude.
\end{enumerate}
\end{proof}

The function $f$ in equation (\ref{eqpoles}) extends to a regular function in a neighborhood of $x+y\s$ if and only if $m \geq 0$. In this case, we can extend the definitions of spherical and isolated multiplicity to transcendental functions. If $m<0$, then we can make use of theorem \ref{transfactorization} to give an analogous definition for the poles of $f$.

\begin{defi}
Let $f, m, n, p_1$ be as in theorem \ref{transfactorization}. If $m\geq0$ then we say that $f$ has \textnormal{spherical multiplicity} $2m$ at $x+y\s$. If, on the contrary, $m<0$ then we say that $f$ has \textnormal{spherical order} $-2m$ at $x+y\s$. In both cases, if $n>0$ then we say that $f$ has \textnormal{isolated multiplicity} $n$ at $p_1$.
\end{defi}

We extend definition \ref{multiplicity} to semiregular functions as:

\begin{defi}
Let $f$ be a semiregular function on $B=B(0,R)$ and let $p \in B$. The \textnormal{(classical) multiplicity} of $f$ at $p$, denoted $m_f(p)$, is the largest $k \in \nn$ such that $f(q) = (q-p)^{*k}*g(q)$ for some semiregular function $g$ on $B$ with $ord_g(p) = 0$, if such a $k$ exists. If not, then $m_f(p) = 0$.
\end{defi}

Remark \ref{relationmult} immediately extends to all semiregular functions $f$ on a ball $B(0,R)$ and all points $p$ in the domain where $f$ is regular.
Similarly, in the case of poles the spherical order and isolated multiplicity of $f$ are related to the (classical) order $ord_f$ as follows.

\begin{rmk}
Let $f$ be a semiregular function on $B(0,R)$ which is not regular at $p = x+yI \in B(0,R)$. Then $f$ has spherical order $2 \max \{ord_f(p),ord_f(\bar p) \}$ at $x+y\s$. If moreover $ord_f(p)>ord_f(\bar p)$, then $f$ has isolated multiplicity $n \geq ord_f(p)-ord_f(\bar p)$ at $\bar p$.
\end{rmk}


\section*{Acknowledgements}

The author acknowledges support by the Unione Matematica Italiana, the MIT-Italy Program and the Department of Mathematics of the Massachusetts Institute of Technology while writing this paper.



\begin{thebibliography}{20}

\bibitem{cohn} P. M. Cohn, \textit{Skew fields. Theory of general division rings}. Encyclopedia of Mathematics and its applications, 57. Cambridge University Press, Cambridge, 1995. 14--18.

\bibitem{cullen} C. G. Cullen, An integral theorem for analytic intrinsic functions on quaternions. \textit{Duke Math. J.} {\bf 32} (1965), 139--148.

\bibitem{openarxiv} G. Gentili, C. Stoppato, The open mapping theorem for quaternionic regular functions, E-print. arXiv:0802.3861v1 [math.CV]

\bibitem{open} G. Gentili, C. Stoppato, The open mapping theorem for quaternionic regular functions, Preprint. Dipartimento di Matematica ``U. Dini'', Universit\`{a} di Firenze, n. 2 (2008).

\bibitem{zeros} G. Gentili, C. Stoppato, Zeros of regular functions and polynomials of a quaternionic variable, Preprint. Dipartimento di Matematica ``U. Dini'', Universit\`{a} di Firenze, n. 1 (2007).

\bibitem{cras} G. Gentili, D. C. Struppa, A new approach to Cullen-regular functions of a quaternionic variable. \textit{C. R. Math. Acad. Sci. Paris} {\bf 342} (2006), 741--744.

\bibitem{advances}  G. Gentili, D. C. Struppa, A new theory of regular functions of a quaternionic variable. \textit{Adv. Math.} {\bf 216} (2007), 279--301.

\bibitem{multiplicity} G. Gentili, D. C. Struppa, On the multiplicity of the zeroes of polynomials with quaternionic coefficients, Preprint. Dipartimento di Matematica ``U. Dini'', Universit\`{a} di Firenze, n. 11 (2007).

\bibitem{lam} T. Y. Lam, \textit{A first course in noncommutative rings}. Graduate Texts in Mathematics, 123. Springer-Verlag, New York, 1991. 261--263.

\bibitem{lam2} T. Y. Lam, \textit{Lectures on modules and rings}. Graduate Texts in Mathematics, 189. Springer-Verlag, New York, 1999. 299--303.

\bibitem{rowen} L. H. Rowen, \textit{Ring theory. Student edition}. Academic press, San Diego, 1991. 272--279.

\bibitem{shapiro} A. Pogorui, M. V. Shapiro, On the structure of the set of zeros of quaternionic polynomials.  \textit{Complex Variables Theory Appl.} {\bf 49} (2004), no. 6, 379--389.

\bibitem{serodio} R. Ser\^{o}dio, L. S. Siu, Zeros of quaternion polynomials. \textit{Appl. Math. Letters} {\bf 14} (2001), 237--239.

\end{thebibliography}
\end{document}